\documentclass{amsart}
\parskip=\smallskipamount

\newtheorem{theorem}{Theorem}[section]
\newtheorem{lemma}[theorem]{Lemma}
\newtheorem{corollary}[theorem]{Corollary}
\newtheorem{proposition}[theorem]{Proposition}

\theoremstyle{definition}
\newtheorem{definition}[theorem]{Definition}
\newtheorem{example}[theorem]{Example}

\newtheorem{problem}[theorem]{Problem}
\newtheorem{remark}[theorem]{Remark}
\newtheorem{question}[theorem]{Question}

\newcommand{\C}{\mathbb{C}}

\newcommand{\N}{\mathbb{N}}

\renewcommand{\P}{\mathbb{P}}
\newcommand{\R}{\mathbb{R}}
\newcommand{\T}{\mathbb{T}}

\newcommand{\cF}{\mathcal{F}}

\newcommand{\cO}{\mathcal{O}}

\def\dibar{\overline\partial}
\def\bs{\backslash}

\numberwithin{equation}{section}

\begin{document}
\title[Runge approximation on convex sets implies the Oka property]
{Runge approximation on convex sets implies \\ the Oka property}
\author{Franc Forstneri\v c}
\address{Institute of Mathematics, Physics and Mechanics, 
University of Ljubljana, Jadranska 19, 1000 Ljubljana, Slovenia}
\email{franc.forstneric@fmf.uni-lj.si}
\thanks{Research supported by grants P1-0291
and J1-6173, Republic of Slovenia.}

%
%
\date{May 5, 2005} 

\begin{abstract}
We prove that the classical Oka property 
of a complex manifold $Y\!$, concerning the existence
and homotopy classification of holomorphic mappings 
from Stein manifolds to $Y\!$,
is equivalent to a Runge approximation property for holomorphic 
maps from compact convex sets in Euclidean spaces to $Y$.
\end{abstract}

\maketitle

\section*{INTRODUCTION} 
Motivated by the seminal works of Oka \cite{O} and Grauert 
(\cite{G1}, \cite{G2}, \cite{G3})
we say that a complex manifold $Y$ enjoys the {\em Oka property} 
if for every Stein manifold $X$, every compact
$\mathcal{O}(X)$-convex subset $K$ of $X$ and every continuous map
$f_0\colon X\to Y$ which is holomorphic in an open neighborhood
of $K$ there exists a homotopy of continuous maps 
$f_t \colon X\to Y$ $(t\in [0,1])$ such that for every $t\in [0,1]$
the map $f_t$ is holomorphic in a neighborhood of $K$
and uniformly close to $f_0$ on $K$, and the map 
$f_1\colon X\to Y$ is holomorphic. 

The Oka property and its generalizations play a central role in 
analytic and geometric problems on Stein manifolds
and the ensuing results are commonly referred to as 
{\em the Oka principle}. Applications include the homotopy classification 
of holomorphic fiber bundles with complex homogeneous fibers 
(the {\em Oka-Grauert principle} \cite{G3}, \cite{Ca}, \cite{HL}) 
and optimal immersion and embedding theorems for Stein manifolds
\cite{EG}, \cite{Sch}; for futher references see 
the surveys \cite{F4} and \cite{Le}.  

In this paper we show that the Oka property 
is equivalent to a Runge-type approximation property
for holomorphic mappings from Euclidean spaces.

%
%
%
%
\begin{theorem} 
\label{Tmain}
If $Y$ is a complex manifold such that any holomorphic map from a 
neighborhood of a compact convex set 
$K\subset \C^n$ $(n\in \N)$
to $Y$ can be approximated uniformly on $K$ 
by entire maps $\C^n\to Y$ 
then $Y$ satisfies the Oka property.
\end{theorem}

The hypothesis in theorem \ref{Tmain} will be referred 
to as the {\em convex approximation property} 
(CAP) of the manifold $Y$. The converse implication is obvious 
and hence the two properties are equivalent: 
$$ 
	{\rm CAP} \Longleftrightarrow {\rm the\ Oka\ property}. 
$$
For a more precise result see theorem \ref{T1.2} below.
An analogous equivalence holds in the parametric case
(theorem \ref{T5.1}), and CAP itself implies the one-parametric 
Oka propery (theorem \ref{T5.3}).

To our knowledge, CAP is the first known characterization of 
the Oka property which is stated purely in terms of holomorphic maps 
from Euclidean spaces and which does not involve additional parameters. 
The equivalence in theorem \ref{Tmain} seems rather striking 
since linear convexity is not a biholomorphically invariant 
property and it rarely suffices to fully describe global complex 
analytic phenomenona. (For the role of convexity in complex analysis 
see  H\"ormander's monograph \cite{H}.)

In the sequel \cite{F8} to this paper it is shown that
CAP of a complex manifold $Y$ also implies the universal 
extendibility of holomorphic maps from closed 
complex submanifolds of Stein manifolds to $Y$
(the {\em Oka property with interpolation}).

We actually show that a rather special class of compact 
convex sets suffices to test the Oka property (theorem 1.2). 
This enables effective applications of the rich theory 
of holomorphic automorphisms of Euclidean spaces 
developed in 1990's, beginning with the works of Anders\'en 
and Lempert \cite{A}, \cite{AL}, thus yielding 
a new proof of the Oka property in several 
cases where the earlier proof relied on sprays 
introduced by Gromov \cite{Gr2}; an example 
are complements of thin (of codimension at 
least two) algebraic subvarieties in certain
algebraic manifolds (corollary \ref{C1.3}).

Theorem \ref{Tmain} partly answers a question, raised by  
Gromov (\cite{Gr2}, p.\ 881, 3.4.(D)), whether 
Runge approximation on a certain class of compact
sets in Euclidean spaces, for example the balls,
suffices to infer the Oka property. 
While it may conceivably be possible to reduce the
testing family to balls by more careful 
geometric considerations, 
we feel that this would not substantially simplify 
the verification of CAP in concrete examples. 

CAP has an essential advantage over the other known sufficient 
conditions when considering unramified holomorphic fibrations
$\pi\colon Y\to Y'$. While it is a difficult problem
to transfer a spray on $Y'$ to one on $Y$ and vice versa, 
lifting an individual map $K\to Y'$ from a convex 
(hence contractible) set $K\subset \C^n$ to a map $K\to Y$ 
is much easier --- all one needs is the Serre fibration
property of $\pi$ and some analytic flexibility condition
for the fibers (in order to find a holomorphic lifting). 
In such case the total space $Y$ satisfies 
the Oka property if and only if the base space $Y'$ does;
this holds in particular if $\pi$ is a holomorphic fiber 
bundle whose fiber satisfies CAP
(theorems \ref{T1.4} and \ref{T1.8}). 
This shows the Oka property for Hopf manifolds, 
Hirzebruch surfaces, complements 
of finite sets in complex tori of dimension $>1$,
unramified elliptic fibrations, etc.

The main conditions on a complex manifold which are known 
to imply the Oka property are complex homogeneity 
(Grauert \cite{G1}, \cite{G2}, \cite{G3}), 
the existence of a dominating spray (Gromov \cite{Gr2}), 
and the existence of a finite 
dominating family of sprays \cite{F2} (def.\ \ref{D2} below).
It is not difficult to see that each of them  
implies CAP --- one uses the given condition 
to linearize the approximation problem and thereby 
reduce it to the classical Oka-Weil approximation theorem 
for section of holomorphic vector bundles over Stein manifolds.
(See also \cite{FP1} and \cite{FP3}. An analogous result 
for algebraic maps has recently been proved in \S 3 of \cite{F7}.) 
The gap between these sufficient conditions and 
the Oka property is not fully understood; 
see \S 3 of \cite{Gr2} and the papers 
\cite{F7}, \cite{F8}, \cite{La}, \cite{La2}.

Our proof of the implication 
CAP$\Rightarrow$Oka property (\S 3 below)
is a synthesis of recent developments 
from \cite{F5} and \cite{F6} where similar
methods have been employed in the construction
of holomorphic submersions. In a typical inductive step
we use CAP to approximate a family of 
holomorphic maps $A\to Y$ from a compact 
strongly pseudoconvex domain $A\subset X$,
where the parameter of the family belongs to 
$\C^p$ $(p=\dim Y)$, by another family of maps 
from a convex bump $B\subset X$ attached to $A$. 
The two families are patched together into a family of 
holomorphic maps $A\cup B \to Y$ by applying a generalized 
Cartan lemma proved in \cite{F5}
(lemma \ref{L2.1} below); this does not require any special 
property of $Y$ since the problem is 
transferred to the source Stein manifold $X$. 
Another essential tool from \cite{F5} allows us 
to pass a critical level of a strongly plurisubharmonic 
Morse exhaustion function on $X$ by reducing the 
problem to the noncritical case for another strongly 
plurisubharmonic function. The crucial part of extending 
a partial holomorphic solution to an attached handle 
(which describes the topological change at a Morse 
critical point)  does not use any condition on 
$Y$ thanks to a Mergelyan-type approximation theorem 
from \cite{F6}.

\section{The main results}
Let $z=(z_1,\ldots,z_n)$ be the coordinates on $\C^n$, with 
$z_j=x_j+iy_j$. Set
$$
	P =\{z\in \C^n\colon |x_j| \le 1,\ |y_j|\le 1,\ j=1,\ldots,n\}.  
								\eqno(1.1)
$$
A {\em special convex set} in $\C^n$ is a compact convex subset of the form
$$ 
	Q=\{z\in P \colon y_n \le h(z_1,\ldots,z_{n-1},x_n)\},   \eqno(1.2)
$$
where $h$ is a smooth (weakly) concave function with values in $(-1,1)$.

We say that a map is holomorphic on a compact set
$K$ in a complex manifold $X$ if it is holomorphic in 
an unspecified open neighborhood of $K$ in $X$; for a homotopy 
of maps the neighborhood should not depend on the parameter.

\begin{definition}  
\label{D1}
A complex manifold $Y$ satisfies 
the {\em $n$-dimensional convex approximation property}
{\rm (CAP$_n$)} if any holomorphic map $f\colon Q\to Y$ 
on a special convex set $Q \subset \C^n$ (1.2) can be 
approximated uniformly on $Q$ by holomorphic maps 
$P\to Y$. $Y$ satisfies $\text{\rm CAP}=\text{\rm CAP}_\infty$ 
if it satisfies ${\rm CAP}_n$ for all $n\in \N$.
\end{definition}

Let $\mathcal{O}(X)$ denote the algebra of all 
holomorphic functions on $X$. A compact set $K$ in $X$ 
is $\mathcal{O}(X)$-convex if for every $p\in X\bs K$
there exists $f\in\cO(X)$ such that 
$|f(p)|> \sup_{x\in K} |f(x)|$. 

\begin{theorem} 
\label{T1.2}
{\bf (The main theorem)}
If $Y$ is a $p$-dimensional complex manifold satisfying 
${\rm CAP}_{n+p}$ for some $n\in \N$ then $Y$ 
enjoys the Oka property  for maps $X\to Y$ from any Stein
manifold with $\dim X\le n$. Furthermore, 
sections $X\to E$ of any holomorphic fiber bundle 
$E\to X$ with such fiber $Y$ satisfy the {\em Oka principle}:
Every continuous section $f_0\colon X\to E$ is homotopic to
a holomorphic section $f_1\colon X\to E$ through a homotopy
of continuous sections $f_t\colon X\to E$ $(t\in [0,1])$; 
if in addition $f_0$ is holomorphic on a compact $\mathcal{O}(X)$-convex 
subset $K \subset X$ then the homotopy $\{f_t\}_{t\in [0,1]}$ 
can be chosen holomorphic and uniformly close to $f_0$ on $K$.
\end{theorem}

Note that the Oka property of $Y$ is just the Oka principle 
for sections of the trivial (product) bundle $X\times Y\to X$
over any Stein manifold $X$. 

We have an obvious implication 
${\rm CAP}_n \Longrightarrow {\rm CAP}_k$ 
when $n>k$ (every compact convex set in $\C^k$ is also
such in $\C^n$ via the inclusion 
$\C^k\hookrightarrow \C^n$), but the converse fails in general
for $n\le \dim Y$ (example \ref{E1}).
An induction over an increasing sequence of cubes 
exhausting $\C^n$ shows that CAP$_n$ is equivalent 
to the Runge approximation of holomorphic maps
$Q\to Y$ on special convex sets (1.2) by entire maps  
$\C^n\to Y$ (compare with the definition of CAP
in the introduction).  

We now verify CAP in  several specific examples. The following 
was first proved in \cite{Gr2} and \cite{F2} by finding a 
dominating family of sprays (see def.\ \ref{D2} below).

\begin{corollary}
\label{C1.3}
Let $p>1$ and let $Y'$ be one of the manifolds $\C^p$, $\C\P_p$ 
or a complex Grassmanian of dimension $p$. If $A \subset Y'$ is a closed
algebraic subvariety of complex codimension at least two then 
$Y=Y'\backslash A$ satisfies the Oka property.
\end{corollary}

\begin{proof}
Let $f\colon Q\to Y$ be a holomorphic map from a special 
convex set $Q\subset P\subset \C^n$ (1.2). 
An elementary argument shows that $f$ can be 
approximated uniformly on $Q$ by algebraic
maps $f' \colon \C^n\to Y'$ such that 
${f'}^{-1}(A)$ is an algebraic subvariety
of codimension at least two which is disjoint from $Q$. 
(If $Y'=\C^p$ we may take a suitable generic polynomial approximation
of $f$, and the other cases easily reduce to this one by the
arguments in \cite{F6}.)  By lemma 3.4 in \cite{F5} there is a 
holomorphic automorphism $\psi$ of $\C^n$ which approximates 
the identity map uniformly in a neighborhood of $Q$ and 
satisfies $\psi(P) \cap {f'}^{-1}(A) =\emptyset$.
The holomorphic map $g = f' \circ \psi \colon \C^n\to Y'$ 
then takes $P$ to $Y= Y'\backslash A$ and it 
approximates $f$ uniformly on $Q$. This proves that 
$Y$ enjoys CAP and hence (by theorem \ref{T1.2}) the Oka property. 
\end{proof}

By methods in \cite{F7} (especially corollary 2.4 and proposition 5.4)
one can extend corollary \ref{C1.3} to any algebraic manifold 
$Y'$ which is a finite union of Zariski open sets 
biregularly equivalent to $\C^p$. Every such manifold
satisfies an approximation property analogous to CAP 
for regular algebraic maps (corollary 1.2 in \cite{F7}).

We now consider unramified holomorphic fibrations,
beginning with a result which is easy to state
(compare with Gromov \cite{Gr2}, 3.3.C' and 3.5.B'',
and L\'arusson \cite{La}, \cite{La2}); 
the proof is given in \S 4.

\begin{theorem}
\label{T1.4}
If $\pi \colon Y\to Y'$ is a holomorphic fiber bundle
whose fiber satisfies {\rm {CAP}} then $Y$ enjoys  
the Oka property if and only if $Y'$ does. This holds in particular
if $\pi$ is a covering projection, or if the fiber of 
$\pi$ is complex homogeneous.
\end{theorem}

\begin{corollary}
\label{C1.5}
Each of the following manifolds enjoys the Oka property:
\begin{itemize}
\item[(i)]   A Hopf manifold.
\item[(ii)]  The complement of a finite set in 
a complex torus of dimension $>1$.  
\item[(iii)] A Hirzebruch surface.
\end{itemize}
\end{corollary}

\begin{proof} (i):
A $p$-dimensional Hopf manifold is a 
holomorphic quotient of $\C^p\backslash \{0\}$ by an 
infinite cyclic group of dilations of $\C^p$ (\cite{BH}, p.\ 225);
since $\C^p\bs \{0\}$ satisfies CAP by corollary \ref{C1.3},
the conclusion follows from theorem \ref{T1.4}.
Note that Hopf manifolds are non-algebraic and even 
non-K\"ahlerian.

(ii):   Every $p$-dimensional torus is a quotient 
$\T^p = \C^p/\Gamma$ where $\Gamma\subset \C^p$ 
is a lattice of maximal real rank $2p$.
Choose finitely many points 
$t_1,\ldots,t_m\in \T^p$ and preimages $z_j\in \C^p$ 
with $\pi(z_j)=t_j$ $(j=1,\ldots,m)$. The discrete set 
$\Gamma' = \cup_{j=1}^m(\Gamma + z_j) \subset \C^p$ 
is tame according to proposition 4.1 in \cite{BL}.
(The cited proposition is stated for $p=2$, but the 
proof remains valid also for $p>2$.) Hence the complement 
$Y= \C^p\backslash \Gamma'$ admits a dominating spray 
and therefore satisfies the Oka property \cite{Gr2}, \cite{FP1}. 
Since $\pi|_Y \colon Y\to \T^p\backslash \{t_1,\ldots,t_m\}$
is a holomorphic covering projection, theorem \ref{T1.4} implies 
that the latter set also enjoys the Oka property. 

The same argument applies if the lattice $\Gamma$ has 
less than maximal rank.

(iii): A {\em Hirzebruch surface} $H_l$ $(l=0,1,2,\ldots)$
is the total space $Y$ of a holomorphic fiber bundle
$Y\to \P_1$ with fiber $\P_1$ (\cite{BH}, p.\ 191;
every Hirzebruch surface is birationally equivalent 
to $\P_2$). Since the base and the fiber are complex homogeneous, 
the conclusion follows from theorem \ref{T1.4}.
\end{proof}

In this paper, an {\em unramified holomorphic fibration} 
will mean a surjective holomorphic submersion 
$\pi\colon Y\to Y'$ which is also a {\em Serre fibration} 
(i.e., it satisfies the homotopy lifting property;
see \cite{W}, p.\ 8). The latter condition holds  
if $\pi$ is a topological fiber bundle in which the 
holomorphic type of the fiber may depend on the base point. 
(Ramified fibrations, or fibrations with multiple fibers, 
do not seem amenable to our methods and will not be discussed;
see example 6.3 and problem 6.7 in \cite{F7}.)
In order to generalize theorem \ref{T1.4}
to such fibration we must assume that the 
fibers of $\pi$ over small open subsets of the base manifold 
$Y'$ satisfy certain condition, analogous to CAP, which allows
holomorphic approximation of local sections. The most general 
known sufficient condition is {\em subellipticity}
\cite{F2}, a generalization of Gromov's ellipticity 
\cite{Gr2}. We recall the relevant definitions.

%
%
%
%
Let $\pi \colon Y\to Y'$ be a holomorphic submersion onto $Y'$. 
For each $y\in Y$ let $VT_y Y =\ker d\pi_y \subset T_y Y$ 
(the {\em vertical tangent space} of $Y$ with respect to $\pi$).
A {\em fiber-spray} associated to $\pi \colon Y\to Y'$ 
is a triple $(E,p,s)$ consisting of a holomorphic vector bundle 
$p\colon E\to Y$ and a holomorphic spray map 
$s\colon E\to Y$ such that for each $y\in Y$ we have $s(0_y)=y$ 
and $s(E_y) \subset Y_{\pi(y)}= \pi^{-1}(\pi(y))$.
A spray on a complex manifold $Y$ is a fiber-spray associated 
to the trivial submersion $Y\to point$.  

%
%
\begin{definition}
\label{D2} 
{\rm (\cite{F2}, p.\ 529)}
A holomorphic submersion $\pi \colon Y\to Y'$
is {\em subelliptic} if each point in $Y'$ has an open neighborhood 
$U\subset Y'$ such that the restricted submersion
$h\colon Y|_U=h^{-1}(U) \to U$ admits finitely many 
fiber-sprays $(E_j,p_j,s_j)$ $(j=1,\ldots,k)$ satisfying  
the domination condition
$$ 
	(ds_1)_{0_y}(E_{1,y}) + (ds_2)_{0_y}(E_{2,y})\cdots 
                     + (ds_k)_{0_y}(E_{k,y})= VT_y Y   \eqno(1.3)
$$
for each $y\in Y|_U$; such a collection of sprays is said to 
be {\em fiber-dominating}. The submersion is {\em elliptic} 
if the above holds with $k=1$. A complex manifold $Y$ is 
(sub-)elliptic if the trivial submersion $Y\to point$ is such. 
\end{definition}

A holomorphic fiber bundle $Y\to Y'$ is (sub-)elliptic when its fiber is such.

\begin{definition}
\label{D3}
A holomorphic map $\pi\colon Y\to Y'$ is a 
{\em subelliptic Serre fibration} if it is a
surjective subelliptic submersion and a Serre fibration.
\end{definition}

The following result is proved in \S 4 below
(see also \cite{La2}).

\begin{theorem} 
\label{T1.8}
If  $\pi \colon Y\to Y'$ is a subelliptic Serre fibration 
then $Y$ satisfies the Oka property if and only if $Y'$ does. 
This holds in particular if $\pi$ is an unramified 
elliptic fibration (i.e., every fiber $\pi^{-1}(y')$ 
is an elliptic curve). 
\end{theorem}

{\em Organization of the paper.}
In \S 2 we state a generalized Cartan lemma 
used in the proof of theorem \ref{T1.2},
indicating how it follows from 
theorem 4.1 in \cite{F5}. Theorem \ref{T1.2} 
(which includes theorem \ref{Tmain}) is proved
in \S 3. In \S 4 we prove theorems \ref{T1.4}
and \ref{T1.8}. In \S 5 we discuss the parametric case and 
prove that CAP implies the one-parametric Oka property
(theorem \ref{T5.3}). \S 6 contains a discussion 
and a list of open problems.

%
%
%
%
\section{A Cartan type splitting lemma}
Let $A$ and $B$ be compact sets in a complex manifold $X$ 
satisfying the following:  
\begin{itemize}
\item[(i)]  $A\cup B$ admits a basis of Stein neighborhoods in $X$, and
\item[(ii)] $\overline {A\backslash B} \cap \overline {B\backslash A}=\emptyset$ 
(the separation property). 
\end{itemize}

Such $(A,B)$ will be called a {\em Cartan pair} in $X$. 
(The definition of a Cartan pair often includes an additional 
Runge condition; this will not be necessary here.)
Set $C=A\cap B$. Let $D$ be a compact set 
with a basis of open Stein neighborhoods in a complex manifold $T$. 
With these assumptions we have the following.

%
%
\begin{lemma}
\label{L2.1}
Let $\gamma(x,t)=(x,c(x,t)) \in X\times T$ $(x\in X,\ t\in T)$
be an injective holomorphic map in an open neighborhood 
$\Omega_C \subset X\times T$ of $C\times D$. If $\gamma$ 
is sufficiently uniformly close to the identity map 
on $\Omega_C$ then there exist open neighborhoods 
$\Omega_A,\Omega_B\subset X\times T$ of $A\times D$, 
respectively of $B\times D$, and injective holomorphic maps 
$\alpha \colon \Omega_A \to X\times T$,
$\beta \colon \Omega_B\to X\times T$ of the
form $\alpha(x,t)=(x,a(x,t))$, $\beta(x,t)=(x,b(x,t))$,
which are uniformly close to the identity map on their 
respective domains and satisfy 
$$   
	\gamma = \beta \circ \alpha^{-1}    
$$
in a neighborhood of $C\times D$ in $X\times T$.
\end{lemma}

In the proof of theorem \ref{T1.2} (\S 3) we shall 
use lemma 2.1 with $D$ a cube in $T=\C^p$
for various values of $p\in \N$. Lemma 2.1 generalizes 
the classical Cartan lemma (see e.g.\ \cite{GR}, p.\ 199) 
in which $A$, $B$ and $C=A\cap B$ are cubes in $\C^n$ 
and $a,b,c$ are invertible  linear functions of $t\in \C^p$ 
depending holomorphically on the base variable.

\begin{proof}
Lemma \ref{L2.1} is a special case of theorem 4.1 in \cite{F5}.
In that theorem we consider a Cartan pair $(A,B)$ in a 
complex manifold $X$ and a nonsingular holomorphic foliation
$\cF$ in an open neighborhood of $A\cup B$ in $X$.
Let $U\subset X$ be an open neighborhood of $C=A\cap B$ in $X$. 
By theorem 4.1 in \cite{F5}, 
every injective holomorphic map $\gamma\colon U\to X$
which is sufficiently uniformly close to the identity map
on $U$ admits a splitting $\gamma = \beta \circ \alpha^{-1}$
on a smaller open neighborhood of $C$ in $X$, 
where $\alpha$ (resp.\ $\beta$) is an injective 
holomorphic map on a neighborhood of $A$ (resp.\ $B$),
with values in $X$. If in addition $\gamma$ preserves
the plaques of $\mathcal{F}$  in a certain finite system 
of foliation charts covering $U$ (i.e., $x$ and $\gamma(x)$ 
belong to the same plaque) then $\alpha$ and $\beta$ can be chosen 
to satisfy the same property. 

Lemma \ref{L2.1} follows by applying this result to the 
Cartan pair $(A\times D, B\times D)$ in $X\times T$,
with $\cF$ the trivial (product) foliation of 
$X\times T$ with leaves $\{x\}\times T$.
\end{proof}

Certain generalizations of lemma \ref{L2.1} are possible
(see \cite{F5}). First of all, the analogous result holds
in the parametric case. 
Secondly, if $\Sigma$ is a closed complex subvariety 
of $X\times T$ which does not intersect $C\times D$ then 
$\alpha$ and $\beta$ can be chosen tangent to 
the identity map to a given finite order along $\Sigma$.
Thirdly, shrinking of the domain is necessary 
only in the directions of the leaves of $\cF$;
an analogue of lemma \ref{L2.1} can be proved for maps 
which are holomorphic in the interior of the respective set
$A,B$, or $C$ and of a H\"older class $\C^{k,\epsilon}$
up to the boundary. (The $\dibar$-problem which 
arises in the linearization is well behaved 
on these spaces.)  We do not state or prove this
generalization formally since it will not 
be needed in the present paper.

%
%
%
%
\section{Proof of theorem \ref{T1.2}}
The proof relies on Grauert's bumping method 
which has been introduced
to the Oka-Grauert problem by Henkin and Leiterer 
\cite{HL} (their paper is based on a preprint 
from 1986), with several additions from 
\cite{F5} and \cite{F6}.

Assume that $Y$ is a complex manifold satisfying CAP.
Let $X$ be a Stein manifold, $K\subset X$ a compact 
$\mathcal{O}(X)$-convex 
subset of $X$ and $f\colon X\to Y$ a continuous map 
which is holomorphic in an open set $U\subset X$ containing $K$.
We shall modify $f$ in a countable sequence of steps to obtain a 
holomorphic map $X\to Y$ which is homotopic to $f$ and approximates 
$f$ uniformly on $K$. (In fact, the entire homotopy will remain
holomorphic and uniformly close to $f$ on $K$.)
The goal of every step is to enlarge the domain of holomorphicity 
and thus obtain a sequence of maps $X\to Y$ which converges 
uniformly on compacts in $X$ to a solution of the problem.

Choose a smooth strongly plurisubharmonic Morse exhaustion 
function $\rho\colon X\to \R$ such that 
$\rho|_K<0$ and $\{\rho\le 0\} \subset U$.
Set $X_c=\{\rho\le c\}$ for $c\in\R$. 
It suffices to prove that 
for any pair of numbers $0\le c_0<c_1$ 
such that $c_0$ and $c_1$ are
regular values of $\rho$, a continuous map 
$f\colon X\to Y$ which is holomorphic on 
(an open neighborhood of) $X_{c_0}$ can 
be deformed by a homotopy of maps 
$f_t\colon X\to Y$ $(t\in [0,1])$ 
to a map $f_1$ which is holomorphic on 
$X_{c_1}$; in addition we require that $f_t$ is 
holomorphic and uniformly as close as required  
to $f=f_0$ on $X_{c_0}$ for every $t\in [0,1]$. 
The solution is then obtained by an obvious induction
as in \cite{FP1}.

There are two main cases to consider:

\smallskip
{\em The noncritical case:} $d\rho\ne 0$ on the set 
$\{x\in X\colon c_0\le \rho(x)\le c_1\}$. 

\smallskip
{\em The critical case:} there is a point $p\in X$
with $c_0<\rho(p)< c_1$ such that $d\rho_p=0$. 
(We may assume that there is a unique such $p$.)
\smallskip

A reduction of the critical case to the noncritical one
has been explained in \S 6 of \cite{F6}, based 
on a technique developed in the construction
of holomorphic submersions of Stein manifolds
to Euclidean spaces \cite{F5}. It is accomplished 
in the following three steps, the first two of which 
do not require any special properties of $Y$.

\smallskip
{\em Step 1:} Let $f\colon X\to Y$ be a continuous map
which is holomorphic in a neighborhood of $X_c=\{\rho\le c\}$ 
for some $c<\rho(p)$ close to $\rho(p)$. 
By a small modification we make $f$ smooth on a totally 
real handle $E$ attached to $X_c$ and passing through the critical point $p$. 
(In suitable local holomorphic coordinates on $X$ near $p$, this handle 
is just the stable manifold of $p$ for the gradient flow of $\rho$,
and its dimension equals the Morse index of $\rho$ at $p$.)
\smallskip

{\em Step 2:} We approximate $f$ uniformly on $X_c\cup E$ 
by a map which is holomorphic in an open neighborhood
of this set (theorem 3.2 in \cite{F6}).
\smallskip

{\em Step 3:} We approximate the map in Step 2 
by a map holomorphic on $X_{c'}$ for some $c'>\rho(p)$.
This extension across the critical level $\{\rho = \rho(p)\}$ 
is obtained by applying the noncritical case for another 
strongly plurisubharmonic function constructed 
especially for this purpose. 
\smallskip

After reaching $X_{c'}$ for some $c'>\rho(p)$ we revert back 
to $\rho$ and continue (by the noncritical case) to the next 
critical level of $\rho$, thus completing the induction step. 
The details can be found in \S 6 in \cite{F5} and \cite{F6}.

It remains to explain the noncritical case; here our proof
differs from the earlier proofs (see e.g.\ \cite{FP1} and \cite{F2}). 
  
Let $z=(z_1,\ldots,z_n)$, $z_j=u_j+iv_j$,
denote the coordinates on $\C^n$,  $n=\dim X$. 
Let $P$ denote the open cube
$$
	P = \{z\in \C^n\colon |u_j|<1,\ |v_j|<1,\ j=1,\ldots,n\}
						\eqno(3.1)
$$
and $P'=\{z\in P\colon v_n=0\}$. Let $A$ be a compact 
strongly pseudoconvex domain with smooth boundary in $X$. 
We say that a compact subset $B\subset X$ is a {\em convex bump} 
on $A$ if there exist an open neighborhood $V\subset X$ of $B$, 
a biholomorphic map $\phi\colon V\to P$ onto the set (3.1) 
and smooth strongly concave functions 
$h, \widetilde   h \colon P'\to [-s,s]$ for some $0<s<1$ such that 
$h\le \widetilde   h$, $h=\widetilde   h$ near 
the boundary of $P'$, and
\begin{eqnarray*}
         \phi(A\cap V) &= \{z\in P \colon  
  	 v_n \le h(z_1,\ldots,z_{n-1},u_n)\}, \cr
             \phi((A\cup B)\cap V) &= \{z\in P \colon  
  	 v_n \le \widetilde   h (z_1,\ldots,z_{n-1},u_n)\}. 
\end{eqnarray*}
We also require that
\begin{itemize}
\item[(i)]
$\overline{A\backslash B}\,\cap\, \overline{B\backslash A}=\emptyset$
(the separation condition), and  
\item[(ii)]
$C=A\cap B$ is Runge in $A$, in the sense that every holomorphic function
in a neighborhood of $C$ can be approximated uniformly on $C$ 
by functions holomorphic in a neighborhood of $A$.
\end{itemize}

\begin{proposition} 
\label{P3.1}
Assume that $A, B\subset X$ 
are as above. Let $Y$ be a $p$-dimensional complex manifold satisfying 
${\rm CAP}_{n+p}$. Choose a distance function $d$ on $Y$ induced by a 
Riemannian metric. For every holomorphic map $f_0 \colon A\to Y$ 
and every $\epsilon >0$ there is a holomorphic map $f_1\colon A\cup B \to Y$ 
satisfying $\sup_{x\in A} d(f_0(x), f_1(x)) <\epsilon $.
The analogous result holds for sections of a holomorphic 
fiber bundle $Z\to X$ with fiber $Y$ which is trivial 
over the set $V\supset B$.
\end{proposition}

If $f_0$ and $f_1$ are sufficiently uniformly close on $A$, 
there clearly exists a holomorphic homotopy from $f_0$ to $f_1$ on $A$. 
If $Y$ satisfies ${\rm CAP}_{N}$ with $N=p+[ \frac{1}{2}(3n+1)]$ 
then we may omit the hypothesis that $C$ be Runge in $A$
(remark \ref{R3.3}). 

Assuming proposition \ref{P3.1} we can complete 
the proof of the noncritical case 
(and hence of theorem \ref{T1.2}) as follows.
By Narasimhan's lemma on local convexi\-fi\-ca\-tion 
of strongly pseudoconvex domains one can obtain a finite sequence 
$X_{c_0}=A_0\subset A_1\subset \ldots\subset A_{k_0}= X_{c_1}$ 
of compact strongly pseudoconvex domains in $X$ such that for every 
$k=0,1,\ldots,k_0-1$ we have $A_{k+1}=A_k\cup B_k$ 
where $B_k$ is a convex bump on $A_k$ (lemma 12.3 in \cite{HL2}).
Each of the sets $B_k$ may be chosen sufficiently small so 
that it is contained in an element of a given open covering 
of $X$. The separation condition (i) for the pair $(A_k,B_k)$, 
introduced just before proposition \ref{P3.1}, is trivial to 
satisfy while (ii) is only a small addition 
(one can use a local convexification of a strongly pseudoconvex 
domain $A$ given by holomorphic functions defined in a 
neighborhood of $A$; see \cite{Fn}, p.\ 530, proposition 1, 
or \cite{HC}, proposition 14). It remains to apply proposition 
\ref{P3.1} inductively  to every pair $(A_k,B_k)$, $k=0,1,\ldots,k_0-1$.
A more detailed exposition of this construction can be
found in \cite{FP1} and \cite{F6}.

This completes the proof of theorem \ref{T1.2} provided 
that proposition \ref{P3.1} holds.

\smallskip
\noindent \em Proof of proposition \ref{P3.1}. \rm 
Choose a pair of numbers $r,r'$, with $0< r'<r<1$, such that 
$\phi(B) \subset r'\overline P$. The set 
$$
	Q := \phi(A\cap V)\cap r\overline P = 
	\{z\in r\overline P\colon  v_n \le h(z',u_n) \}
$$	 
is a special convex set in $\C^n$ (1.2) with respect to 
the closed cube $r\overline P \subset\C^n$, and  
the set $C=A\cap B$ is contained in 
$Q_0 := \phi^{-1}(Q)\subset X$. 

By the hypothesis $f_0$ is holomorphic in an open neighborhood $U\subset X$ 
of $A$. Set $F_0(x)=(x,f_0(x)) \in X\times Y$ for $x\in U$.

\begin{lemma}
\label{L3.2}
There are a neighborhood $U_1\subset U$ of $A$ in $X$, a 
neighborhood $W\subset \C^p$ of $0\in \C^p$ and a 
holomorphic map $F(x,t)=(x,f(x,t)) \in X\times Y$,
defined for $x\in U_1$ and $t\in W$, such that 
$f(\cdotp,0)=f_0$ and $f(x,\cdotp) \colon W\to Y$ 
is injective holomorphic for every $x$ in a neighborhood 
of $C=A\cap B$. 
\end{lemma}

\begin{proof}
The set $F_0(U)$ is a closed Stein submanifold of 
the complex manifold $U\times Y$ and hence it admits 
an open Stein neighborhood in $U\times Y$ according
to \cite{Siu}. Let $\pi_X\colon X\times Y\to X$ denote the 
projection $(x,y)\to x$. 
The set $E =\ker d\pi_X$  is a holomorphic vector 
subbundle of rank $p=\dim Y$ in the tangent 
bundle $T(X\times Y)$, consisting of all vectors
$\xi\in T(X\times Y)$ which are tangent to the fibers of $\pi_X$. 

Since the set $Q_0$ is contractible, the bundle $E$ is 
trivial over a neighborhood of $F_0(Q_0)$ in $X\times Y$ 
and hence is generated there by $p$ holomorphic sections, 
i.e., vector fields tangent to the fibers of $\pi_X$. 
Since $C$ is Runge in $A$, these sections 
can be approximated uniformly on $F_0(C)$ by 
holomorphic sections $\xi_1,\ldots,\xi_p$
of $E$, defined in a neighborhood of $F_0(A)$ in $X\times Y$,
which still generate $E$ over a neighborhood of $F_0(C)$. 
The flow $\theta^j_t$ of $\xi_j$ is well defined
for sufficiently small $t\in \C$. The map 
$$
	F(x,t_1,\ldots,t_p)= \theta^1_{t_1}\circ\cdots \circ 
	\theta^p_{t_p}\circ F_0(x)  \in X\times Y,
$$
defined and holomorphic for $x$ in a neighborhood of $A$
and for $t=(t_1,\ldots,t_p)$ in a neighborhood of the origin
in $\C^p$, satisfies lemma 3.2. 
\end{proof}

\begin{remark} 
\label{R3.3}
The restriction of a rank $p$ 
holomorphic vector bundle $E$ to an $n$-dimensional Stein 
manifold is generated by $p+[\frac{1}{2}(n+1)]$ sections 
(lemma 5 in \cite{Fo}, p.\ 178). Without assuming
that $C$ is Runge in $A$ this gives a proof of 
lemma \ref{L3.2} if $Y$ satisfies ${\rm CAP}_{N}$ with
$N=p+[\frac{1}{2}(3n+1)]$.
\end{remark}

We continue with the proof of proposition 3.1.
Let $F$ and $W$ be as in lemma 3.2.
Choose a closed cube $D$ in $\C^p$ centered at $0$, 
with $D \subset W$. The set 
$\widetilde   Q :=Q \times D \subset \C^{n+p}$ is a special
convex set of the form (1.2) with respect to the 
closed cube $\widetilde P := r\overline P\times D\subset \C^{n+p}$,
and the map $\widetilde f(z,t):=f(\phi^{-1}(z),t) \in Y$ 
is holomorphic in a neighborhood of $\widetilde Q$. 

Since $Y$ is assumed to satisfy CAP$_{n+p}$, 
we can approximate $\widetilde f$ uniformly on a neighborhood 
of $\widetilde Q$ by entire maps $\widetilde  g\colon \C^{n+p}\to Y$. 
(This is the only place in the proof where CAP is used.) 
The holomorphic map 
$$
	g(x,t) := \widetilde g(\phi(x),t)\in Y, \qquad x\in V,\ t\in \C^p
$$  
then approximates $f$ uniformly in a 
neighborhood of $Q_0 \times D$ in $X\times\C^p$.
Since $f(x,\cdotp)  \colon W\to Y$ is injective holomorphic 
for every $x$ in a neighborhood of $C$ (lemma \ref{L3.2}),
choosing $g$ to approximate $f$ sufficiently well we obtain
a (unique) injective holomorphic map 
$\gamma(x,t)=(x,c(x,t)) \in X\times \C^p$, defined 
and uniformly close to the identity map in an open 
neighborhood $\Omega \subset X\times \C^p$ of $C \times D$, 
such that    
$$   
	f(x,t) = (g\circ\gamma)(x,t) = g(x,c(x,t)), \quad (x,t)\in \Omega. 
								\eqno(3.2)
$$
If the approximation of $f$ by $g$ is sufficiently close then 
$\gamma$ is so close to the identity map that we can apply
lemma 2.1 to obtain a decomposition $\gamma=\beta\circ\alpha^{-1}$,
with $\alpha(x,t)=(x,a(x,t))$, $\beta(x,t)=(x,b(x,t))$
holomorphic and close to the identity maps in 
their respective domains $\Omega_A \supset A\times D$, 
$\Omega_B\supset B\times D$. From (3.2) we obtain 
$$  
	f(x,a(x,t)) = g(x,b(x,t)), \quad  (x,t)\in C\times D.
$$
Setting $t=0$, the two sides define a holomorphic map 
$f_1\colon A\cup B \to Y$ which approximates $f_0=f(\cdotp,0)$ 
uniformly on $A$ (since $a(x,0)\approx 0$ for $x\in A$). 

This proves proposition 3.1 for maps $X\to Y$.  
The very same proof applies to sections of a holomorphic 
fiber bundle $Z\to X$ with fiber $Y$ which is trivial over 
the set $V\supset B$; this is no restriction since all 
convex bumps in the inductive construction can be chosen 
small enough to insure this condition.

%
%
%
%
\section{Proof of theorems \ref{T1.4} and \ref{T1.8}} 
We begin by proving theorem \ref{T1.8}.
Let $\pi\colon Y\to Y'$ be a subelliptic Serre fibration
(definition \ref{D3}). Assume first that $Y'$ satisfies CAP.  
Let $f\colon U\to Y$ be a holomorphic map from 
an open convex subset $U\subset \C^n$.
Let $K\subset L$ be compact convex sets in $U$, 
with $K\subset {\rm Int}\, L$. 
Set $g=\pi \circ f\colon U\to Y'$. 

Since $Y'$ satisfies CAP, there is an entire map $g_1\colon\C^n\to Y'$ 
which approximates $g$ uniformly on $L$. 

By lemma 3.4 in \cite{F6} there exists for every $x\in U$ 
a holomorphic retraction $\rho_x$ of an open neighborhood 
of the fiber $R_x=\pi^{-1}(g_1(x)) \subset Y$ in the manifold $Y$ 
onto $R_x$, with $\rho_x$ depending holomorphically 
on $x\in U$. If $g_1$ is sufficiently uniformly close 
to $g$ on $L$ then for every $x\in L$ the point $f(x)$ 
belongs to the domain of $\rho_x$, and hence we can define 
$f_1(x)= \rho_x(f(x))$ for all $x\in L$. The map $f_1$ is 
then holomorphic on a neighborhood of $K$ in $X$, 
it approximates $f$ uniformly on $K$, and it satisfies 
$\pi\circ f_1=g_1$ (i.e., $f_1$ is a lifting of $g_1$). 

Since $\pi\colon Y\to Y'$ is a Serre fibration and the set 
$K \subset \C^n$ is convex, $f_1$ extends to a continuous map 
$f_1 \colon \C^n\to Y$ which is holomorphic in a neighborhood 
of $K$ and satisfies $\pi\circ f_1= g_1$ on $\C^n$
(hence $f_1$ is a global lifting of $g_1$). 

Since $g_1$ is holomorphic and $\pi$ is a subelliptic submersion, 
theorem 1.3 in \cite{F3} shows that we can homotopically 
deform $f_1$ (through liftings of $g_1$)
to a global holomorphic lifting $\widetilde   f \colon \C^n\to Y$ 
of $g_1$ (i.e., $\pi\circ \widetilde f = g_1$) such that 
$\widetilde   f|_K$ approximates $f_1|_K$, and hence $f|_K$.
(In our case $\pi$ is unramified and the quoted theorem 
from \cite{F3} is an immediate consequence of theorem 1.5 in \cite{FP2}.) 
This shows that $Y$ satisfies CAP and hence the Oka property. 

Conversely, assume that $Y$ satisfies CAP. 
Choose a holomorphic map $g\colon K\to Y'$ from a compact 
convex set $K\subset \C^n$.  Since $\pi$ is a Serre fibration
and $K$ is contractible, there is a continuous lifting 
$f_0\colon K \to Y$ with $\pi\circ f_0=g$. Since $\pi$ 
is a subelliptic submersion, theorem 1.3 in \cite{F3} gives a 
homotopy of liftings $f_t\colon K\to Y$ $(t\in [0,1])$, 
with $\pi\circ f_t=g$ for every $t\in [0,1]$, such that 
$f_1$ is holomorphic on $K$. 

By CAP of $Y$ we can approximate $f_1$ uniformly on $K$ 
by entire maps $\widetilde   f\colon \C^n\to Y$. 
The map $\widetilde  g := \pi\circ \widetilde   f \colon \C^n\to Y'$ 
is then entire and it approximates $g$ uniformly on $K$. 
Thus $Y'$ satisfies CAP. 

Note that contractibility of $K$ was essential in the 
last part of the proof. 

Every unramified elliptic fibration $\pi\colon Y\to Y'$
without exceptional (and multiple) fibers is elliptic
in the sense of Gromov \cite{Gr2} (definition 
\ref{D2} above). Indeed, every fiber 
$Y_y=\pi^{-1}(y)$ $(y\in Y')$ is an elliptic curve,
$Y_y=\C/\Gamma_y$, and the lattice $\Gamma_y \subset \C$
is defined over every sufficiently small open subset 
$U\subset Y'$ by a pair of generators $a(y)$, $b(y)$ 
depending holomorphically on $y$. A dominating 
fiber-spray on $Y|_U$ is obtained by pushing down to
$Y|_U$ the $\Gamma_y$-equivariant spray on $U\times \C$
defined by $((y,t),t') \in U\times \C\times \C \to (y,t+t') \in U\times \C$. 

The proof of theorem \ref{T1.4} follows the same
scheme; in this case we do not need to refer to 
\cite{FP2} but can instead use theorem \ref{T1.2} 
in this paper.

%
%
%
%
\section{The parametric convex approximation property} 
We recall the notion of the {\em parametric Oka property}
(POP) which is essentially the same as 
Gromov's ${\rm Ell}_\infty$ property (\cite{Gr2}, \S 3.1; 
see also theorem 1.5 in \cite{FP2}).

Let $P$ be a compact Hausdorff space 
(the parameter space) and $P_0$ a closed subset of $P$ 
(possibly empty)
which is a strong deformation retract of some neighborhood in $P$. 
In applications $P$ is usually a polyhedron and $P_0$ a subpolyhedron. 

Given a Stein manifold $X$ and a compact $\mathcal{O}(X)$-convex 
subset $K$ in $X$, we consider a continuous map $f\colon X\times P\to Y$ 
such that for every $p\in P$ the map $f^p=f(\cdotp,p)\colon X\to Y$ 
is holomorphic in an open neighborhood of $K$ in $X$ (independent of $p\in P$), 
and for every $p\in P_0$ the map $f^p$ is holomorphic on $X$. 
We say that $Y$ satisfies the {\em parametric Oka property} 
(POP) if for every such data $(X,K,P,P_0,f)$ there is a homotopy 
$f_t\colon X\times P\to Y$ $(t\in [0,1])$, consisting
of maps satisfying the same properties as $f_0=f$, such that 
\begin{itemize}
\item[---] the homotopy is fixed on $P_0$ 
(i.e., $f_t^p=f^p$ when $p\in P_0$ and $t\in [0,1]$), 
\item[---] $f_t$ approximates $f$ uniformly on 
$K\times P$ for all $t\in [0,1]$, and
\item[---]  $f_1^p \colon X\to Y$ is holomorphic  
for every $p\in P$.
\end{itemize}
Recall that POP is implied by ellipticity
\cite{Gr2}, \cite{FP1} and subellipticity \cite{F3}. 

We say that a complex manifold $Y$ satisfies the 
{\em parametric convex approximation property}
(PCAP) if the above holds for every special convex set 
$K$ of the form (1.2) in $X=\C^n$ for any $n\in \N$. 

\begin{theorem} 
\label{T5.1} 
If a complex manifold $Y$ satisfies {\rm PCAP} then 
it also satisfies the parametric Oka property
(and hence ${\rm PCAP}\, \Longleftrightarrow\, {\rm POP}$).
\end{theorem}

Theorem \ref{T5.1} is obtained by following the proof of 
theorem \ref{T1.2} (\S 3) but using the requisite tools 
with continuous dependence on the parameter $p\in P$. 
Precise arguments of this kind can be found in \cite{FP1}, \cite{FP2} 
and we leave out the details. 
For an additional equivalence involving interpolation
conditions see theorem 6.1 in \cite{F8}.

An analogue of theorem \ref{T1.8} holds for ascending/descending 
of the parametric Oka property (POP) in a subelliptic 
Serre fibration $\pi\colon Y\to Y'$. The implication 
$$ 
	{\rm POP\ of\ }Y' \,\Longrightarrow\,  {\rm POP\ of\ }Y 
$$
holds for any compact Hausdorff parameter space $P$ 
and is proved as before by using the parametric versions 
of the relevant tools.  However, we can prove the converse implication 
only for a {\em contractible} parameter space $P$,
the reason being that we must lift a map $K\times P\to Y'$ 
(with $K$ a compact convex set in $\C^n$)
to a map $K\times P\to Y$. (See also corollary 6.2 in \cite{F8}.)

\begin{question}
\label{CAPimpliesPCAP}
To what extent does CAP imply PCAP?
\end{question}

We indicate how CAP$\Rightarrow$PCAP can be proved 
for sufficiently simple parameter spaces. For simplicity 
let $P$ be a closed cube in $\R^k$ and $P_0=\emptyset$, 
although the argument applies in more general situations.
We identify $\R^k$ with $\R^k\times\{i0\}^k \subset\C^k$.

Let $K\subset \C^n$ be a special compact convex
set (def.\ \ref{D1}), $U\subset \C^n$ an open neighborhood 
of $K$, and $f\colon U\times P \to Y$ a 
continuous map such that $f_p=f(\cdotp,p)$ is holomorphic 
on $U$ for every fixed $p\in P$. By the assumed CAP property
of $Y$ we can approximate $f_p$ for every fixed $p\in P$
uniformly on $K$ by a map with values in $Y$ 
which is holomorphic in an open neighborhood
of $K\times \{p\}$ in $\C^n\times \C^k$.
Patching these holomorphic approximations 
by a smooth partition of unity in the $p$-variable 
we approximate the initial map $f$ by another one,
still denoted $f$, which is smooth in all variables 
and is holomorphic in the $x$ variable for every fixed $p\in P$.

The graph of $f$ over $U \times P$ is a smooth CR submanifold 
of $\C^{n+k}\times Y$ foliated by $n$-dimensional complex manifolds, 
namely the graphs of $f_p \colon U\to Y$ for $p\in P$. 
By methods similar to those in \cite{FCL} it can be seen 
that the graph of $f$ over $K\times P$ admits an 
open Stein neighborhood $\Omega$ in $\C^{n+k}\times Y$.
Embedding $\Omega$ into a Euclidean space $\C^N$ and 
applying standard approximation methods for CR functions
(and a holomorphic retraction of a tube around 
the submanifold $\Omega\subset \C^N$ onto $\Omega$) 
we can approximate $f$ as closely as desired on 
$K\times P$ by a holomorphic map $\widetilde f$, 
defined in an open neighborhood of 
$K\times P$ in $\C^n \times \C^p$.

The cube $P \subset \R^k\subset \C^k$ 
admits a basis of cubic neighborhoods in $\C^k$.
(By a `cube' we mean a Cartesian product of intervals 
in the coordinate axes.) 
The product of $K$ with a closed cube in $\C^k$
is a special compact convex set in $\C^{n+k}$.
Applying the CAP property of $Y$ to the map 
$\widetilde f$ we see that $Y$ satisfies PCAP for 
the parameter space $P$. 

If $P =[0,1]\subset \R$ and 
the maps $f_0=f(\cdotp,0)$ and $f_1=f(\cdotp, 1)$ 
(corresponding to the endpoints of $P$) are holomorphic, 
the above construction can be performed so that these
two maps remain unchanged, thereby showing that 
the basic CAP implies the one-parametric CAP. 
Joined with theorem \ref{T5.1} this gives

\begin{theorem}
\label{T5.3}
If a complex manifold $Y$ enjoys {\rm CAP} then a 
homotopy of maps $f_t \colon X\to Y$ $(t\in [0,1])$ 
from a Stein manifold $X$ for which $f_0$ and $f_1$ are 
holomorphic can be deformed with fixed ends to 
a homotopy consisting of holomorphic maps.
\end{theorem}

Theorem \ref{T5.3} also follows from theorem 1.1 in \cite{F8} 
to the effect that {\em CAP implies the Oka property with interpolation}. 
Indeed, extending the homotopy $f_t \colon X\to Y$ $(t\in [0,1])$  to all 
values $t\in \C$ by precomposing with a retraction $\C\to [0,1] \subset \C$ 
we obtain a continuous map $F\colon X\times \C \to Y$,
$F(x,t)=f_t(x)$, whose restriction to the complex submanifold 
$X_0= X\times \{0,1\}$ of $X\times \C$ is holomorphic.
By \cite{F8} there is a homotopy $F_s \colon X\times\C \to Y$
($s\in [0,1]$), with $F_0=F$, which remains fixed on $X_0$
and such that $F_1$ is holomorphic. The restriction
of $F_1$ to $X\times [0,1]$ is a homotopy 
from $f_0$ to $f_1$ consisting of holomorphic
maps $F_1(\cdotp,t)$ $(t\in [0,1])$.

%
%
%
%
\section{Discussion, examples and problems}
It was pointed out by Gromov \cite{Gr2}
that the existence of a dominating spray on a 
complex manifold $Y$ is a precise way of saying
that $Y$ admits  many holomorphic 
maps from Euclidean spaces; since every Stein
manifold $X$ embeds into a Euclidean space, this
also implies the existence of many holomorphic maps 
$X\to Y$ and hence it is natural to expect
that $Y$ enjoys the Oka property (and it does). 

The same philosophy 
justifies CAP which is another way of asserting
the existence of many holomorphic maps $\C^N\to Y$. 
Indeed, CAP is the restriction of the Oka property 
(which refers to maps from {\em any} Stein manifold $X$
to $Y$, with uniform approximation on {\em any} holomorphically 
convex subset $K$ of $X$) to model pairs
--- the special compact convex sets in $X=\C^n$. 
For a discussion of this {\em localization principle} 
see remark 1.10 in \cite{F7}.

CAP is in a precise sense 
opposite to the hyperbolicity properties
expressed by nonvanishing of Kobayashi-Eisenman metrics.
More precisely, ${\rm CAP}_1$ is an opposite property to 
{\em Kobayashi-Brody hyperbolicity} \cite{Kb}, \cite{Br}
which excludes nonconstant entire  maps $\C\to Y$; 
more generally, ${\rm CAP}_n$ for $n\le \dim Y$ is opposite 
to the $n$-dimensional  measure hyperbolicity 
\cite{E}. The property ${\rm CAP}_p$  with 
$p=\dim Y$ implies the existence of 
{\em dominating} holomorphic maps $\C^p\to Y$; 
if such $Y$ is compact, it is not of Kodaira 
general type \cite{CG}, \cite{Kd}, \cite{KO}. 
For a further discussion see \cite{F7}.

The property ${\rm CAP}_n$ for $n\ge \dim Y$ 
is also reminiscent of the Property ${\mathrm S}_{\mathrm n}$, 
introduced in \cite{F6}, which requires that 
any holomorphic submersion $f\colon K\to Y$ 
from a special compact convex set $K\subset \C^n$
is approximable by entire submersions $\C^n\to Y$. 
By theorem 2.1 in \cite{F6}, Property ${\mathrm S}_{\mathrm n}$ 
of $Y$ implies that holomorphic submersions from any $n$-dimensional 
Stein manifold to $Y$ satisfy the homotopy principle,
analogous to the one which was proved for smooth submersions
by Gromov \cite{Gr1} and Phillips \cite{P}. 
The similarity is not merely apparent --- our proof 
of theorem \ref{T1.2}
in this paper conceptually unifies the construction of 
holomorphic maps with the construction of 
holomorphic submersions in \cite{F5} and \cite{F6}.

\begin{example}
\label{E1}
For every $1\le k\le p$ there exists a $p$-dimensional 
complex manifold which satisfies ${\rm CAP}_{k-1}$ but not
${\rm CAP}_{k}$. 
\end{example}

Indeed, for $k=p$ we can take $Y=\C^p\backslash A$
where $A$ is a discrete subset of $\C^p$ which is
{\em rigid} in the sense of Rosay and Rudin 
(\cite{RR}, p.\ 60), i.e., every holomorphic
map $\C^p \to \C^p$ with maximal rank $p$ at some
point intersects $A$ at infinitely many points.
Thus ${\rm CAP}_p$ fails, but ${\rm CAP}_{p-1}$ holds
since a generic holomorphic map $\C^{p-1}\to\C^p$ 
avoids $A$ by dimension reasons. 
For $k < p$ we take $Y=\C^p\backslash \phi(\C^{p-k})$ 
where $\phi\colon \C^{p-k} \hookrightarrow \C^p$ is a proper holomorphic 
embedding whose complement is $k$-hyperbolic (every entire map 
$\C^k\to \C^p$ whose range omits $\phi(\C^{p-k})$ has rank $<k$; 
such maps were constructed in \cite{F1}); 
again CAP$_k$ fails but CAP$_{k-1}$ holds by dimension reasons. 
Another example is $Y= (\C^{k}\backslash A) \times \C^{p-k}$ 
where $A$ is a rigid discrete set in $\C^{k}$.

\smallskip
We conclude by mentioning a few open problems.

\begin{problem}
Do the CAP$_n$ properties stabilize at some integer, i.e.,
is there a $p\in \N$ depending on $Y$ (or perhaps only
on $\dim Y$) such that 
${\rm CAP}_p \,\Longrightarrow\, {\rm CAP}_n$
for all $n> p$~? Does this hold for $p=\dim Y$~?
\end{problem}

\begin{problem}
\label{Q1} 
Let $B$ be a closed ball in $\C^p$ for some $p\ge 2$.
Does $\C^p\backslash B$ satisfy CAP (and hence the Oka property)~?
Does $\C^p\backslash B$ admit any nontrivial sprays~? 
\end{problem}

The same problem makes sense for every compact convex set 
$B\subset \C^p$. What makes this problem particularly intriguing is the
absence of any obvious obstruction; indeed, $\C^p\backslash B$ is 
a union of Fatou-Bieberbach domains \cite{RR}.

\begin{problem}
\label{Q2}
{\rm (Gromov \cite{Gr2}, p.\ 881, 3.4.(D))}
Suppose that every holomorphic map from a ball $B\subset \C^n$
to $Y$ (for any $n\in \N$) can be approximated by entire maps 
$\C^n\to Y$. Does $Y$ enjoy the Oka property~? 
\end{problem}

\begin{problem}
\label{Q3}
Let $\pi\colon Y\to Y_0$ be a holomorphic
fiber bundle. Does the Oka property of $Y$ imply 
the Oka property of the base $Y_0$ and of the fiber~?
\end{problem}

\medskip
{\em Acknowledgement.}
Some of the ideas and techniques in this paper originate in my 
joint works with Jasna Prezelj whom I wish to thank for her indirect 
contribution. I also thank Finnur L\'arusson and Edgar Lee Stout 
for helpful discussions, and the referee for thoughtful 
remarks which helped me to improve the presentation.

\bibliographystyle{amsplain}

\end{document}